\documentclass{amsart}

\newtheorem{theo}{Theorem}
\newtheorem{corol}{Corollary}

\providecommand{\nequiv}{\not\equiv}

\def\COne{C_{+1}}
\def\CTwo{C_{-1}}

\begin{document}

\title{On prime chains}

\author{Douglas S. Stones}
\address{School of Mathematical Sciences, Monash University, VIC 3800 Australia}
\email{douglas.stones@sci.monash.edu.au}

\begin{abstract}
Let $b$ be an odd integer such that $b \equiv \pm 1 \pmod 8$ and let $q$ be a prime with primitive root $2$ such that $q$ does not divide $b$.  We show that if $(p_k)_{k=0}^{q-2}$ is a sequence of odd primes such that $p_k=2p_{k-1}+b$ for all $1 \leq k \leq q-2$, then either (a)~$q$ divides $p_0+b$, (b)~$p_0=q$ or (c)~$p_1=q$.
\end{abstract}

\maketitle

For integers $a,b$ with $a \geq 1$, a sequence of primes $(p_k)_{k=0}^{\lambda-1}$ such that $p_k=ap_{k-1}+b$ for all $1 \leq k \leq \lambda-1$ is called a \emph{prime chain} of length $\lambda$ based on the pair $(a,b)$.  This follows the terminology of Lehmer \cite{Lehmer1965}.  The value of $p_k$ is given by
\begin{equation}\label{GeneralEq}
p_k=a^k p_0+b\frac{(a^k-1)}{(a-1)}
\end{equation}
for all $0 \leq k \leq \lambda-1$.

For prime chains based on the pair $(2,+1)$, Cunningham \cite[p.~241]{Cunningham1907} listed three prime chains of length $6$ and identified some congruences satisfied by the primes within prime chains of length at least $4$.  Prime chains based on the pair $(2,+1)$ are now called \emph{Cunningham chains of the first kind}, which we will call $\COne$ chains, for short.  Prime chains based on the pair $(2,-1)$ are called \emph{Cunningham chains of the second kind}, which we will call $\CTwo$ chains.

We begin with the following theorem which has ramifications on the maximum length of a prime chain; it is a simple corollary of Fermat's Little Theorem.  A proof is also given by L\"{o}h \cite{Loh1989}.  Moser \cite{Moser1952} once posed Theorem~\ref{THBasicStruct}, with $a,b,p_0 \geq 1$, as an exercise, for which he received fourteen supposedly correct proofs.

\begin{theo}\label{THBasicStruct}
Let $(p_k)_{k \geq 0}$ be an infinite sequence for which $p_k=a p_{k-1}+b$ for all $k \geq 1$.  Then the set $\{p_k\}_{k \geq 0}$ is either finite or contains a composite number.
\end{theo}

There are some choices of $(p_0,a,b)$ that are uninteresting.  For example, if $b=-(a-1)p_0$, then the prime chain is $(p_0,p_0,\dots)$.  In fact, if $p_i=p_j$ for any distinct $i,j$ then $\{p_k\}_{k=0}^{\infty}$ will be periodic, with period dividing $|p_i-p_j|$.  Also if $\gcd(a,b)>1$ then the sequence could only possibly be of length $1$, since $\gcd(a,b)$ divides $a p_0+b$.

In this article we will therefore assume that $(p_k)_{k=0}^{\infty}$ is a strictly increasing sequence.  Theorem~\ref{THBasicStruct} implies that no choice of $(p_0,a,b)$ will give rise to a prime chain of infinite length.  However, this raises the question, how long can a prime chain be?  Green and Tao \cite{GreenTao2008} proved that, for all $\lambda \geq 1$, there exists a prime chain of length $\lambda$ based on the pair $(1,b)$ for some $b$.  Lehmer \cite{Lehmer1965} remarked that Dickson's Conjecture \cite{Dickson1904}, should it be true, would imply that there are infinitely many prime chains of length $\lambda$ based on the pair $(a,b)$, with the exception of some inappropriate pairs $(a,b)$.

Discussions about searching for Cunningham chains were given by Lehmer \cite{Lehmer1965}, Guy \cite[Sec. A7]{Guy2004}, Loh \cite{Loh1989} and Forbes \cite{Forbes1999}.  Tables of Cunningham chains are currently being maintained by Wikipedia \cite{WikipediaCunningham} and Caldwell \cite{Caldwell1999}.

In this article, we will frequently deal with primes, denoted either $p$ or $q$, that have a primitive root $a$.  We therefore introduce the following terminology for brevity.  If $a$ is a primitive root modulo $q$ then we will write $a \bigtriangleup q$ and if $q$ is prime and $a \bigtriangleup q$, we will call $q$ an $a \bigtriangleup$-\emph{prime}.

We begin with the following theorem, which slightly improves \cite[Thm~1]{Lehmer1965}.

\begin{theo}\label{Theorem1}
Let $q$ be an $a \bigtriangleup$-prime such that $q$ does not divide $b$.  Suppose $(p_k)_{k=0}^{q-2}$ is a prime chain based on the pair $(a,b)$.  Then $q$ divides $p_0(a-1)+b$ or $q=p_k$ for some $0 \leq k \leq q-2$.
\end{theo}

\begin{proof}
To begin, note that $a \nequiv 0,1 \pmod q$ since $a \bigtriangleup q$.  Suppose $q$ does not divide $p_0(a-1)+b$.  Let $S=\{p_k\}_{k=1}^{q-2}$ and let $S_q=\{p_k \pmod q\}_{k=1}^{q-2}$.  If $p_i \equiv p_j \pmod q$ then \[a^i p_0 + b \frac{(a^i-1)}{(a-1)} \equiv a^j p_0 + b \frac{(a^j-1)}{(a-1)} \pmod q\] by \eqref{GeneralEq} and so \[\frac{a^i}{a-1} \big(p_0(a-1)+b\big) \equiv \frac{a^j}{a-1} \big(p_0(a-1)+b\big) \pmod q\] since $a \nequiv 1 \pmod q$.  Since $q$ does not divide $p_0(a-1)+b$ we find $a^i \equiv a^j \pmod q$ implying that $i \equiv j \pmod {q-1}$, since $a \bigtriangleup q$.  Therefore $|S_q| = q-1$.

If $-b/(a-1) \pmod q \in S_q$ then for some $i$, \[a^i p_0 + b \frac{(a^i-1)}{(a-1)} \equiv \frac{-b}{a-1} \pmod q\] by \eqref{GeneralEq}, implying that $p_0 \equiv -b/(a-1) \pmod q$ contradicting our initial assumption.  Hence $S_q=\{0,1,2,\dots,q-1\} \setminus \{-b/(a-1)\}$.  Since $q$ does not divide $b$, we find that $0 \in S_q$ and therefore $q$ divides an element of $S$.  But since $S$ contains only primes, therefore $q \in S$.
\end{proof}

To show that Theorem~\ref{Theorem1} is the ``best possible'' in at least one case, we identify the prime chain $(7,11,23,59,167,491)$ of length $\lambda=q-1=6$ based on $(a,b)=(3,-10)$.  Here $-b/(a-1)=10/2 \equiv 5 \pmod 7$ while $p_0 \equiv 0 \pmod 7$.  This raises the question, when can there exist a prime chain of length $q-1$, for $q$, $a$ and $b$ satisfying the conditions of Theorem~\ref{Theorem1} while $p_0 \nequiv -b/(a-1) \pmod q$?  In the next section, we will find that prime chains of this form, when $b \equiv \pm 1 \pmod 8$, are exceptional, which includes Cunningham chains of both kinds.

% \section{Cunningham chains}\label{SECunChain}

Cunningham \cite[p.~241]{Cunningham1907} claimed that a $\COne$ chain $(p_k)_{k=0}^{\lambda-1}$ of length $\lambda \geq 4$ must have (a)~each $p_k \equiv -1 \pmod 3$ and (b)~each $p_k \equiv -1 \pmod 5$.  However, condition (b) is incorrect for the prime chain $(2,5,11,23,47)$.  In Theorem~\ref{THBpm1Mod8} we will prove that there are no counter-examples to Cunningham's condition (b) when $p_0>5$.  Lehmer \cite{Lehmer1965} stated that $\COne$ chains $(p_i)_{i=0}^{\lambda-1}$ of length $\lambda \geq 10$ have $p_0 \equiv -1 \pmod {2 \cdot 3 \cdot 5 \cdot 11}$.  Loh \cite{Loh1989} showed that $\CTwo$ chains $(p_i)_{i=0}^{\lambda-1}$ of length $\lambda \geq 12$ have $p_0 \equiv 1 \pmod {2 \cdot 3 \cdot 5 \cdot 11 \cdot 13}$.  In Corollary~\ref{COC1C2Chains} we will generalise this list of results to prime chains based on $(2,b)$ for all odd integers $b \equiv \pm 1 \pmod 8$.

For any odd prime $s$ let $o_s(2)$ denote the multiplicative order of $2$ modulo $s$.  Let $\mathbb{N}=\{1,2,\dots\}$.  We make use of the following Legendre symbol identities, which can be found in many elementary number theory texts, for example \cite{LeVeque1996}.  For odd prime $q$
\begin{equation}\label{EQLegendre}
\left(\frac{2}{q}\right)=\begin{cases} 1 & \text{if } q \equiv \pm1 \pmod 8 \\ -1 & \text{if } q \equiv \pm 3 \pmod 8 \end{cases} \hspace{10pt} \text{ and } \hspace{10pt} \left(\frac{a}{q}\right) \equiv a^{\frac{q-1}{2}} \pmod q.
\end{equation}
We are now ready to state and prove the main theorem.

\begin{theo}\label{THBpm1Mod8}
Let $b$ be an odd integer such that $b \equiv \pm 1 \pmod 8$ and let $q$ be a $2 \bigtriangleup$-prime that does not divide $b$.  Suppose $(p_k)_{k=0}^{q-2}$ is a prime chain based on the pair $(2,b)$.  Then either (a)~$q$ divides $p_0+b$, (b)~$p_0=q$, (c)~$p_1=q$ or (d)~$p_0=2$.
\end{theo}

\begin{proof}
Suppose $p_0$ is an odd prime and is of the form $p_0=2m-b$ for some integer $m$.  So $p_k=2^{k+1} m-b$ for all $0 \leq k \leq q-2$ by \eqref{GeneralEq}.  If $q$ does not divide $p_0+b$, then Theorem~\ref{Theorem1} implies that $q=p_k=2^{k+1}m-b$ for some $0 \leq k \leq q-2$.  If $q=2^{k+1} m-b$ where $k \geq 3$, then \[1=\left(\frac{2}{q}\right) \equiv 2^{\frac{q-1}{2}} \pmod q\] by \eqref{EQLegendre} since $b \equiv \pm 1 \pmod 8$.  However, this contradicts that $2 \bigtriangleup q$.  Hence $q=p_0$ or $q=p_1$.
\end{proof}

We can now deduce the following corollary, for which we make use of the fact that contiguous subsequences of prime chains are themselves prime chains to find a large divisor for $p_0-1$.  Let $\lambda \in \mathbb{N}$.

\begin{corol}\label{COC1C2Chains}
Let $(p_k)_{k=0}^{\lambda-1}$ be a prime chain based on $(2,b)$ for an odd integer $b \equiv \pm 1 \pmod 8$.  Let \[Q=\{q \leq \lambda+1: q\text{ is a prime and } 2 \bigtriangleup q\} \cup \{2\} \setminus \{p_0,p_1\}.\]  If $p_0 \geq 3$, then $p_k+b$ is divisible by every $q \in Q$ for all $0 \leq k \leq \lambda-1$.
\end{corol}

\begin{proof}
We know that $2$ divides each $p_k+b$ since both $p_k$ and $b$ are odd.  So let $q \in Q \setminus \{2\}$.  If $p_0 \equiv -b \pmod q$ then $q$ divides $p_k+1$ for all $0 \leq k \leq \lambda-1$ since $2 \times (-b)+b \equiv b \pmod q$.  Observe that $(p_i)_{i=0}^{q-2}$ is a prime chain of length $q-1$ for all $q \in Q$.  The result therefore follows from Theorem~\ref{THBpm1Mod8}.
\end{proof}

The $2\bigtriangleup$-primes are given by Sloane's \cite{Sloane} A001122 as $3,5,11,13,19,29,37,53$, and so on.  It would also be of interest to know if an analogue of Corollary~\ref{COC1C2Chains} holds for other non-trivial values of $(a,b)$.  The techniques in this article use the Legendre symbol identity \eqref{EQLegendre} which requires $a=2$ and $b \equiv \pm 1 \pmod 8$, so they are not easily extended to encompass other pairs $(a,b)$.

Corollary~\ref{COC1C2Chains} does not hold for when $p_0=2$.  For example, $(2,5,11,23,47)$ is a prime chain based on $(2,1)$.  In fact, the subsequences $(2,5,11,23)$ and $(5,11,23,47)$ also illustrate why we need to exclude $p_0$ and $p_1$ from $Q$ in Corollary~\ref{COC1C2Chains}.

Finally, the author would like to thank Hans Lausch for valuable feedback.

\bibliographystyle{amsplain}
% \bibliography{LatinB4}

\providecommand{\bysame}{\leavevmode\hbox to3em{\hrulefill}\thinspace}
\providecommand{\MR}{\relax\ifhmode\unskip\space\fi MR }
\providecommand{\MRhref}[2]{%
  \href{http://www.ams.org/mathscinet-getitem?mr=#1}{#2}
}
\providecommand{\href}[2]{#2}

\end{document}